\documentclass[11pt]{amsart}

\usepackage[latin1]{inputenc}
\usepackage{amssymb,amsmath,amsfonts}
\usepackage{enumerate}

\newtheorem{MainTheorem}{Theorem}

\newtheorem{Proposition}{Proposition}[section]

\newtheorem{Lemma}[Proposition]{Lemma}
\newtheorem{Theorem}[Proposition]{Theorem}

\newtheorem{Remark}[Proposition]{Remark}

\DeclareMathOperator{\vol}{Vol}
\DeclareMathOperator{\kl}{Kl}
\DeclareMathOperator{\Val}{Val}
\DeclareMathOperator{\re}{Re}

\DeclareMathOperator{\Kl}{Kl}
\DeclareMathOperator{\Gr}{Gr}

\DeclareMathOperator{\SL}{SL}
\DeclareMathOperator{\GL}{GL}

\DeclareMathOperator{\SO}{SO}
\DeclareMathOperator{\D}{D}

\newcommand{\R}{\mathbb{R}}
\newcommand{\C}{\mathbb{C}}

\newcommand{\K}{\mathcal{K}}

\newcommand{\func}[5]{\ensuremath{\begin{array}{cccl}
#1:&#2&\longrightarrow&#3\\&#4&\mapsto&#5\end{array}}}

\title{Difference bodies in complex vector spaces}
\author{Judit Abardia} 
\address{Institut f\"ur Mathematik, Goethe-Universit\"at Frankfurt, 
Robert-Mayer-Str. 10, 60054 Frankfurt, Germany}

\email{abardia@math.uni-frankfurt.de}

\begin{document}

\begin{abstract}
A complete classification is obtained of continuous, translation invariant, Minkowski valuations on an $m$-dimensional complex vector space which are covariant under the complex special linear group. 
\end{abstract}
\thanks{Supported by DFG grant BE 2484/3-1
}
\date{\today}

\maketitle 
%-------------------------------------------------------------------------

\section{Introduction}

The classification of real- or body-valued valuations satisfying certain natural properties has attracted a lot of attention in the last years. The first fundamental classification result dates back to 1957, when Hadwiger classified the continuous, translation invariant real-valued valuations which are also invariant under the rotations of the Euclidean space. Since then many generalizations of this result have been obtained. 

We denote by $V$ a real vector space of dimension $n$ and by $\K(V)$ the space of compact convex bodies in $V$. An operator $Z:\K(V)\to(A,+)$ with $(A,+)$ an abelian semi-group is called a {\em valuation} if it satisfies the following additivity property $$Z(K\cup L)+ Z(K\cap L)=Z(K)+Z(L),$$ for all $K,L\in\K(V)$ such that $K\cup L\in\K(V)$.

The classical case consists of taking $(A,+)$ as the real numbers with the usual sum. A particular class of real-valued valuations consists of those which are continuous -- with respect to the Hausdorff topology -- and translation invariant, i.e. $Z(K+x)=Z(K)$ for every $x\in V$. Some of the most important and recent results on the theory of continuous translation invariant valuations can be found in  \cite{alesker_mcullenconj01, alesker_fourier, klain00,ludwig_reitzner99, mcmullen80}. This theory has been extended to the more general framework of manifolds instead of a real vector space, see for instance \cite{alesker04_product,bernig_broecker07}. 
Apart from the continuity and the translation invariance of a real-valued valuation, we can impose invariance under some group acting transitively on the sphere (for instance, the unitary group). Then, we always get a finite dimensional real vector space (see \cite{alesker_survey07}). Its dimension, a basis and the arising integral geometry have been studied intensively. For some references on this direction see \cite{alesker03_un,bernig_g2, bernig_fu_hig,bernig_fu_solanes,fu06}.

Some other important particular cases of valuations are given, for instance, when considering the vector space of symmetric tensors (see \cite{alesker99,alesker_bernig_schuster, hug_schneider_schuster_b,ludwig03, mcmullen97} for more information on tensor-valued valuations), or $(\K(V),+)$ with $+$ the Minkowski sum between two convex bodies (i.e. $K+L=\{x+y\,:\,x\in K,\,y\in L\}$). Convex body valued valuations with addition the Minkowski sum are called {\em Minkowski valuations}.

In this paper, we are interested in dealing with Minkwoski valuations. Some results on Minkowksi valuations not described in this paper can be found, for instance, in \cite{haberl10,haberl11,kiderlen, ludwig06_survey, ludwig10, schuster10, schuster.wannerer12, wannerer12}. Some papers dealing with convex geometry, but working in a complex vector space as ambient space -- as we do -- instead of a real vector space are \cite{koldobsky11,koldobsky_koenig_zymonopoulou08, koldobsky.paouris.zymonopoulu}. 

Two fundamental properties of Minkowski valuations are the contravariance and the covariance with respect to the special linear group $\SL(V,\R)$. A valuation $Z:\K(V)\to\K(V^*)$ is {\em $\SL(V,\R)$-contravariant} if $$Z(gK)=g^{-*}Z(K),\quad\forall g\in\SL(V,\R),$$
where $V^*$ denotes the dual space of $V$ and $g^{-*}$ denotes the inverse of the dual map of $g$. 

A valuation $Z:\K(V)\to\K(V)$ is {\em $\SL(V,\R)$-covariant} if $$Z(gK)=gZ(K),\quad\forall g\in\SL(V,\R).$$

An example of a continuous, translation invariant Minkowski valuations which is $\SL(V,\R)$-contravariant is the projection body operator. For $K\in\K(V)$ the {\em projection body $\Pi K$ of $K$} has support function 
$$h(\Pi K, u)=\frac{n}{2}V(K,\dots,K,[-u,u]),\quad u\in V,$$
where $V(K,\dots,K,[-u,u])$ denotes the mixed volume with $(n-1)$ copies of $K$ and one copy of the segment joining $u$ and $-u$. 

The projection body was introduced in the 19th century by Minkowski and since then it has been widely studied (see, for instance, the books \cite{gardner_book06,koldobsky, leichtweiss_book98,schneider_book93,thompson_book96}). In the framework of the classification results of Minkowski valuations, Ludwig proved in \cite{ludwig02} that the projection body operator is the only (up to a positive constant) continuous Minkowski valuation which is translation invariant and $\SL(V,\R)$-contravariant. In \cite{abardia.bernig} a complex version of this result was shown. The result is as follows
\begin{Theorem}[\cite{abardia.bernig}]\label{thm_contravariant}
Let $W$ be a complex vector space of complex dimension $m \geq 3$. A map $Z:\mathcal{K}(W) \to \mathcal{K}(W^*)$ is a continuous, translation invariant and $\SL(W,\mathbb{C})$-contravariant Minkowski valuation if and only if there exists a convex body $C\subset \C$ such that $Z=\Pi_C$, where $\Pi_C K \in \mathcal{K}(W)$ is the convex body with support function 
\begin{equation}\label{complexProjection}h(\Pi_{C}K,w)=V(K[2m-1],C w),\quad\forall w \in W,\end{equation}
where $C w:=\{cw\,|\,c\in C\subset\C\}$. 
Moreover, $C$ is unique up to translations.
\end{Theorem}

For the covariant case, Ludwig proved in \cite{ludwig_2005} that the difference body is the unique (up to a positive constant) continuous Minkowski valuation which is translation invariant and $\SL(V,\R)$-covariant. In fact, she classified the continuous, $\SL(V,\R)$-covariant Minkowski valuations (not necessarily translation invariant). The {\em difference body of a convex body $K\in\K(V)$} is defined by $$\D\!K=K+(-K),$$ where $-K$ denotes the reflection of $K$ at the origin. 

\medskip
In this paper we study the continuous Minkowksi valuations in a complex vector space $W$ which are translation invariant and $\SL(W,\C)$-covariant. 
Our main result gives a classification of these valuations. 

\begin{MainTheorem}\label{thm_t1}
Let $W$ be a complex vector space of complex dimension $m \geq 3$. A map $Z:\mathcal{K}(W) \to \mathcal{K}(W)$ is a continuous, translation invariant and $\SL(W,\mathbb{C})$-covariant Minkowski valuation if and only if there exists a convex body $C\subset\C$ such that $Z=\D_C$, where
$\D_C\!K \in \mathcal{K}(W)$ is the convex body with support function 
\begin{equation} \label{eq_def_pi_c}
 h(\D_{C}\!K,\xi)=\int_{S^1}h(\alpha K,\xi)dS(C,\alpha),\quad\forall \xi \in W^*,
\end{equation}
where $dS(C,\cdot)$ denotes the area measure of $C$, and $\alpha K=\{\alpha  k\,:\, k\in K\subset W\}$ with $\alpha\in S^1\subset\C$. Moreover, $C$ is unique up to translations.
\end{MainTheorem}

\smallskip
The hypothesis $m\geq 3$ in Theorem \ref{thm_t1} cannot be omitted. In Section \ref{dim2cov} we give for $m=2$ another family of valuations satisfying all the properties and we characterize the continuous, translation invariant Minkowski valuations which are $\SL(W,\C)$-covariant and have fixed degree of homogeneity. We also show that the continuous, translation invariant, $\SL(W,\C)$-contravariant Minkowski valuations with degree of homogeneity 1 are precisely the ones introduced in \cite[Proposition 3.3]{abardia.bernig}.

\subsection*{Acknowledgments}
I would like to express my gratitude to Andreas Bernig for all the enlightening conversations during the preparation of this work and his useful remarks.

%---------------------------------------------------------------------------------------

\section{Background and conventions}

We denote by $V$ a real vector space of dimension $n$ and by $W$ a complex vector space of complex dimension $m$. The space of compact convex bodies in $V$ (resp.\! in $W$) is denoted by $\K(V)$ (resp.\! $\K(W)$). The dual vector space of $V$ (resp.\! $W$) is denoted by $V^*$ (resp.\! $W^*$).

For more information about the notions introduced here we refer to \cite{gardner_book06,gruber_book,schneider_book93}.

\subsection{Support function}

Let $K \in \mathcal{K}(V)$. The {\em support function of $K$} is given by
\begin{align*}
 h_K:V^* & \to \R, \\
\xi & \mapsto \sup_{x \in K}\langle \xi,x\rangle,
\end{align*}
where $\langle \xi,x\rangle$ denotes the pairing of $\xi \in V^*$ and $x\in V$.

The support function is $1$-homogeneous (i.e. $h_K(t\xi)=th_K(\xi)$ for all $t \geq 0$) and subadditive (i.e.
$h_K(\xi+\eta) \leq h_K(\xi)+h_K(\eta)$ for all $\xi,\eta\in V^*$). Moreover, if a function on $V^*$ is $1$-homogeneous and subadditive, then it is the support function of a unique compact convex set $K \in \K(V)$ (cf. \cite[Theorem 1.7.1]{schneider_book93}). 
We also write $h(K,\xi)$ for $h_{K}(\xi)$.

The support function is also linear with respect to the Minkowski sum on $\K(V)$ and has the following important property 
\begin{equation}\label{eq_adjoint}
 h(gK,\xi)=h(K,g^*\xi),\quad \forall \xi \in V^*,\, g \in \GL(V,\R).
\end{equation}
In a complex vector space $W$ this equality holds for $g\in\GL(W,\C)$. In particular, for $\alpha\in\C$ and $K\in\K(W)$ we can interpret $\alpha K=g K$ with $g=\alpha \mathrm{Id}\in\GL(W,\C)$, where $\mathrm{Id}$ denotes the identity matrix. Hence, we have $$h(\alpha K,\xi)=h(K,\alpha^*\xi),$$ 
where $\alpha^*$ denotes $g^*=\overline\alpha\mathrm{Id}$.

The vector space spanned by all support functions has the following density property (cf. \cite[Lemma 1.7.9]{schneider_book93}).
\begin{Lemma}[\cite{schneider_book93}]
\label{differenceSupport}
Every twice-differentiable function on the sphere is the difference of two support functions. 

In particular, the real vector space spanned by the differences of support functions (restricted to $S^{n-1}$) is dense in the space $C(S^{n-1})$ of continuous functions on the sphere (with the maximum norm).
\end{Lemma}

\subsection{Surface area measure and Minkowski's theorem}
Let $K\in\K(V)$, $V$ endowed with a scalar product, and $\omega\subset S^{n-1}$ a Borel subset of $S^{n-1}$. The {\em surface area measure of $K$} is given by $$S(K,\omega)=\vol_{n-1}(\{x\in\partial K\,: \textrm{ an outward unit normal of } x \textrm{ is in } \omega\}).$$

Note that if $K\in\K(V)$ is a polytope, then the surface area measure is a discrete measure: the sum of point masses at the outward unit normal vectors to the facets of $K$, with weight the surface area of the corresponding facet. 

Minkowski's existence theorem gives necessary and sufficient conditions for a positive measure on $S^{n-1}$ to be the surface area measure of some convex body (cf. \cite[Theorem 7.1.2]{schneider_book93}).
\begin{Theorem}[Minkowski's existence theorem]\label{minkowskiExistence}Let $\mu$ be a positive finite Borel measure on $S^{n-1}$. Then, $\mu$ is the surface area measure of some convex body $K\subset V$ with non-empty interior if and only if $\mu$ is not concentrated on any great subsphere of $S^{n-1}$ and 
\begin{equation}\label{zero}\int_{S^{n-1}}u d\mu(u)=0.\end{equation}
\end{Theorem}

\subsection{Translation invariant valuations}
Let $\Val$ denote the Banach space of real-valued, translation invariant, continuous valuations on $V$.

A valuation $\phi \in \Val$ is called {\it homogeneous of degree $k$} if $\phi(tK)=t^k\phi(K)$ for all $t \geq 0$; {\em even} if $\phi(-K)=\phi(K)$ for all
$K$; and {\em odd} if $\phi(-K)=-\phi(K)$. The subspace of even (resp.\! odd) 
valuations of degree $k$ is denoted by $\Val_k^+$ (resp.\! $\Val_{k}^-$). 

\begin{Theorem}[McMullen \cite{mcmullen77}]
\begin{equation} \label{eq_mcmullen_dec}
\Val=\bigoplus_{\substack{k=0,\ldots,n\\ \varepsilon=+,-}} \Val_k^\varepsilon.
\end{equation} 
\end{Theorem}

In \cite{klain95} Klain (see also \cite{klain00}) gives the following description of even translation invariant valuations. For simplicity, we fix a Euclidean scalar product on $V$. Let $\phi \in \Val_k^+$ and let $E$ be a $k$-dimensional subspace of $V$. Klain proved that $\phi|_E$ is a multiple of the volume on $E$, i.e.
\begin{displaymath}
\phi(K)=\kl_\phi(E) \vol(K), \quad \forall K \in \K(E).
\end{displaymath}
The function $\kl_\phi:\Gr_k(V) \to \R$, where $\Gr_k(V)$ the Grassmannian mani\-fold of all $k$-dimensional subspaces in $V$, is called the {\it Klain function of $\phi$}. 

\begin{Theorem}[Klain's injectivity theorem \cite{klain95}] Let $\phi\in\Val_{k}^+$. Then $\phi$ is uniquely determined by its Klain function $\kl_\phi \in C(\Gr_kV)$. 
\end{Theorem}

The group $\GL(V)$ acts naturally on $\Val$ by
\begin{displaymath}
g\mu(K)=\mu(g^{-1}K), \quad g \in \GL(V,\R),\, K \in \mathcal{K}(V). 
\end{displaymath}

A valuation $\mu\in\Val$ is called {\em smooth} if the map $g \mapsto g\mu$ from the Lie group $\GL(V,\R)$ to the Banach space
$\Val$ is smooth. The subspace of smooth valuations is denoted by
$\Val^{sm}$, it is a dense subspace in $\Val$. 
We will use that if $\mu \in \Val_k^{sm,+}$, then the Klain function of $\mu$ is a smooth function on $\Gr_kV$. See \cite{alesker_survey07, alesker_bernig, bernig_broecker07} for more information on smooth valuations.

\subsection{Valuations and distributions}
Let $\mathcal{E}$ denote the space of continuous $1$-homogeneous functions defined on $V^*$. Let $K\subset V^*$ be a compact convex body containing the origin in its interior. Let us endow $\mathcal{E}$ with the supremum norm restricted to $K$ in $V^*$, i.e. $\|f\|_{K}=\sup\{|f(\xi)|\,:\,\xi\in K\}$. Then, for every $K$, $L$ compact convex bodies containing the origin in its interior,  the norms $\|\cdot\|_{K}$, $\|\cdot\|_{L}$ are equivalent in $\mathcal{E}$ and it becomes a Banach space.

Let $\mathcal{D}$ denote the space of the functions in $\mathcal{E}$ which are smooth on $V^*\setminus\{0\}$.

Goodey and Weil \cite{goodey_weil84} give a representation of a continuous, translation invariant, real-valued valuation of homogeneity degree one in terms of a distribution on the sphere $S^{n-1}$. We need the following special case.

\begin{Theorem}[\cite{goodey_weil84}]\label{theoGoodeyWeil} Let $\phi:\K(V)\to\R$ be a continuous, translation invariant valuation which is homogeneous of degree $1$. Then, there exists a unique distribution $T$ on $\mathcal{D}$ which can be extended to the Banach subspace of $\mathcal{E}$ generated by the support functions $h_{K}$ for every $K\in \K(V)$ in such a way that
$$\phi(K)=T(h_{K}).$$
\end{Theorem}

%----------------------------------------------------------------------------------------------------------

\section{Proof of Theorem \ref{thm_t1}}

\begin{Lemma}\label{lemma_homog}
Let $W$ be a complex vector space of complex dimension $m \geq 3$. Let $Z:\K(W) \to \K(W)$ be a continuous, translation invariant, $\SL(W,\mathbb{C})$-covariant Minkowski valuation with degree of homogeneity $k$, $1< k\leq 2m-1$. Then $ZK=\{0\}$, $\forall K \in \K(W)$.
\end{Lemma}

\proof
Let $Z$ be a Minkowski valuation of degree $k$ satisfying the hypothesis of the lemma. Define the operator $\tilde Z:\K(W)\to\K(W)$ by 
\begin{displaymath}
\tilde Z(K):=\int_{S^1} \int_{S^1} q_1 Z(q_2K)dq_1dq_2.
\end{displaymath}
It satisfies $\tilde Z(qK)=\tilde Z(K)$ and $q\tilde Z(K)=\tilde Z(K)$ for all $q\in S^1$ and $K\in\K(W)$. We say that $\tilde Z$ is an {\em $S^1$-bi-invariant} valuation.

$\tilde Z$ inherits all the desired properties from $Z$ and it turns out to be a continuous, translation invariant and $\SL(W,\C)$-covariant Minkowski valuation of degree $k$.
In order to prove the lemma, it suffices to show that there cannot exist a non-trivial $S^1$-bi-invariant valuation satisfying the hypothesis. We denote again this valuation by $Z$. 

Let $g\in \GL(W,\mathbb{C})$ and write $g=g_{0} t q$ with $g_{0}\in\SL(W,\mathbb{C})$, $t \in \mathbb{R}_{>0}, q \in S^1$. Using the $S^1$-bi-invariance and the homogeneity of degree $k$ of $Z$ we have
$$Z(gK)=Z(g_{0}t q K)=t^{k}g_{0}qZ(K)=t^{k-1}gZK,$$ and it follows that
\begin{equation}\label{eq_covariance_gln}Z(gK)=|\!\det g|^{\frac{k-1}{m}}gZK,\quad\forall g\in\GL(W,\mathbb{C}).\end{equation} 

We distinguish two cases.

\noindent {\bf Case $k=m+1$.}
Let $e_1,\ldots,e_m$ be a complex basis of $W$ and $e^1,\ldots,e^m$ its dual basis. We denote by $E$ the $(m+1)$-dimensional real subspace generated by $e_{1},\ldots,e_{m},ie_{1}$.

Let $g \in \GL(W,\mathbb{C})$ be defined by $ge_j=\lambda_j e_j$ with 
$\lambda_1,\ldots,\lambda_m \in \R_{>0}$. Note that $g$ fixes $E$. Let $D=\lambda_1^2 \prod_{j=2}^m \lambda_j$ be the determinant of the restriction of $g$ to $E$ (considered as an element of $\GL(E,\R)$). Let $j\in\{2,\dots,m\}$ and $\xi=e^j$ or $\xi=ie^j$. Using \eqref{eq_covariance_gln} we get 
$$h(ZgK,\xi)=h(ZK,g^{*}\xi)  |\!\det g| =h(ZK,\xi)  \lambda_j |\!\det g|.$$
On the other hand,  by Klain's result, the restriction of $h(Z(\cdot),\xi)$ to
$E$ is a multiple of the $(m+1)$-dimensional volume. Thus, for every $K\in\K(E)$
$$h(ZgK,\xi)=\vol(gK)\Kl(E)=D\vol(K)\Kl(E)=D h(ZK,\xi).$$
Consequently, 
$$h(ZK,\xi)\lambda_{1}\prod_{j=2}^m\lambda_{j}=h(ZK,\xi)\lambda_{j}\prod_{j=2}^m\lambda_{j},\quad\forall\lambda_{1},\dots,\lambda_{m},$$
which implies $$h(ZK,e^j)=h(ZK,ie^j)=0,\quad \forall j\neq 1,\,K\in\K(E).$$
Hence, the support function $h:=h_{ZK}$ vanishes on all lines $\R \cdot e^j$, $\R \cdot ie^j, j=2,\ldots,m$. 
Since $ZK=-ZK$, this implies that $ZK$ is a two-dimensional convex body contained in the space generated by $\{e_{1},ie_{1}\}$.

Let now $g\in\GL(W,\C)$ be defined by $ge_{1}=\alpha e_{1}$, $\alpha=x+iy\in\C$, $ge_{j}=\lambda_{j}e_{j}$, $\lambda_{j}\in\R_{>0}$, $j=\{2,\dots,m\}$. The determinant $D$ of the restriction of $g$ to $E$ is $$D=(x^2+y^2)\lambda_{2}\ldots\lambda_{m},$$
and 
$$|\!\det g|=|\alpha|\lambda_{2}\ldots\lambda_{m}=\sqrt{x^2+y^2}\lambda_{2}\ldots\lambda_{m}.$$
Choosing $\alpha$ with $|\alpha|=1$ we get, for every $K\in\K(E)$
$$h(Z(gK),e^{1})=|\!\det g|h(ZK,\overline\alpha e^{1})=\lambda_{2}\ldots\lambda_{m}h(ZK,\overline\alpha e^1),$$
$$h(Z(gK),e^1)=Dh(ZK,e^1)=\lambda_{2}\ldots\lambda_{m}h(ZK,e^1).$$
Thus, \begin{equation}\label{disc}h(ZK,e^1)=h(ZK,\alpha e^1),\quad\forall\alpha\in S^1,\,K\in\K(E),\end{equation}
and $ZK$ is a disc of radius $r(K)$ contained in the complex line generated by $e_{1}$.

Let $K_{0}\subset E$ be the parallelotope $[0,e_{1}]+[0,ie_{1}]+[0,e_{2}]+\cdots+[0,e_{m}]$ which we denote by $[e_{1},ie_{1},e_{2},\dots,e_{m}]$, and let $K=[w_{1},iw_{1},w_{2}\dots,w_{m}]$ be a parallelotope with $w_{1}=\alpha e_{1}$, $\alpha\in S^1$. We claim that 
\begin{equation}\label{zkk0}h(ZK,e^1)=c|\!\det (w_{1},\dots,w_{m})|,\end{equation}
where $c=h(ZK_{0},e^{1})$. Indeed, using the continuity of both sides of \eqref{zkk0} it is enough to prove it when $w_{1},\dots,w_{m}$ are linearly independent over $\C$. In this case, we can define $g\in\GL(W,\C)$ by $ge_{j}=w_{j}$, $j=1,\dots,m$,  and from \eqref{eq_covariance_gln} and  \eqref{disc} we have
$$h(ZK,e^1)=h(Z(gK_{0}),e^1)=|\!\det g|h(ZK_{0},g^*e^{1})=c|\!\det(w_{1},\dots,w_{m})|.$$

Let us fix a Hermitian scalar product on $W$ such that $e_{1},\dots,e_{m}$ constitutes an orthonormal basis. 

Let $W_{0}$ be the $(m-1)$-dimensional complex subspace of $W$ generated by $\{e_{2},\dots,e_{m}\}$.
Now, let us define a valuation $\phi:\K(W_{0})\to\R$ by
$$\phi(K')=h(Z[e_{1},ie_{1},K'],e^1),$$ where $[e_{1},ie_{1},K']$ denotes the product of the parallelotope $[e_{1},ie_{1}]$ and $K'\subset W_{0}$. Note that both convex sets lie in orthogonal spaces.

Define $H\subset SU(W)$ as the stabilizer of $SU(W)$ at $e_{1}$. We have $H\cong SU(W_{0})\cong SU(m-1)$. If $m\geq 3$, then $H$ acts transitively on the unit sphere of $W_{0}$. 

By \eqref{zkk0}, $\phi$ is $SU(W_{0})$-invariant. Alesker established in \cite[Proposition 2.6]{alesker_survey07} that if $G$ is a compact subgroup of the orthogonal group acting transitively on the unit sphere of a vector space, then each $G$-invariant translation invariant continuous valuation is smooth. Thus, $\phi$ is a smooth valuation. In particular, the Klain function of $\phi$ is a smooth function.

Let us consider the smooth curve $\gamma:\R\to\K(W_{0})$ given by 
$$\gamma(t)=[\cos t e_{2}+\sin t ie_{3},e_{3},\dots,e_{m}].$$
For these convex sets,
\begin{align*}\phi(\gamma(t))&=h(Z[e_{1},ie_{1},\cos t e_{2}+\sin t ie_{3},e_{3},\dots,e_{m}],e^1)\\&=c|\!\det(e_{1},\cos t e_{2}+\sin t ie_{3},e_{3},\dots,e_{m})|=c|\!\cos t|,\end{align*}
which is smooth only if $c=0$.

Hence, we get $h(ZK,e^1)=0$ and from \eqref{disc} we have $h(ZK,\alpha e^1)=0$ for all $K\in\K(E)$ and $\alpha\in S^1$. Thus, $r(K)=0$ (the radius of $ZK$) and by Klain's injectivity theorem we have $Z\equiv \{0\}$.

\noindent {\bf Case $1< k\leq m$ or $m+1 < k\leq 2m-1$.}
The proof of this case is completely analogous to the proof of the contravariant case in \cite[Lemma 3.2]{abardia.bernig} and we do not reproduce it here.
The main idea of the proof was to use the same matrices $g\in\GL(W,\C)$ defined in the previous case. Using \eqref{eq_covariance_gln} and the fact that the power of $|\!\det g|$ is not an integer one obtains that $Z$ must be the trivial valuation. 
\endproof

\begin{Remark}\label{remarkHomog}{\em If $Z:\K(W)\to\K(W)$ is a continuous, translation invariant, $\SL(W,\C)$-covariant  Minkowski valuation of degree $2m$ (resp.\! 0), then the support function of the image is a multiple -- depending on the direction -- of the volume (resp.\! the Euler characteristic) and it can be proved as before that it must be the trivial valuation.}
\end{Remark}

%------------------------

\proof[Proof of Theorem \ref{thm_t1}] 
We assume first that $\D_{C}$ is defined as in \eqref{eq_def_pi_c} and we prove that it satisfies all the stated properties. 

The function on the right hand side of \eqref{eq_def_pi_c} is a support function since $h(\alpha K,\cdot)$ is a support function for every $\alpha$ and $dS(C,\cdot)$ is a positive measure.
Hence $\D_{C}\!K$ is a convex body on $W$ for every $C\in\K(\C)$. 
 
In order to show that $\D_C$ is a Minkowski valuation we use the additivity of the support function in its first argument.  
Let $K,L\in\K(W)$ with $K\cup L\in\K(W)$. Then, $K\cup L+K\cap L=K+L$ (cf. \cite[Lemma 3.1.1]{schneider_book93}) and it follows
\begin{align*}
 h(\D_C(K\cup L) & +\D_C(K\cap L),\xi)  = h(\D_C(K\cup L),\xi)+h(\D_C(K\cap L),\xi)\\
&= \int_{S^1}h(K\cup L+K\cap L,\alpha^*\xi)dS(C,\alpha)\\
&=\int_{S^1}h(K+L,\alpha^*\xi)dS(C,\alpha)\\
& = h(\D_CK,\xi)+h(\D_CL,\xi),
\end{align*}
which implies the valuation property of $\D_C$.

The continuity of $\D_{C}$ follows from the continuity of the support function. 

To prove that $\D_{C}$ is translation invariant we use the only if part of Theorem \ref{minkowskiExistence}. 
Indeed, for $u\in W$ it follows
\begin{align*}
h(\D_{C}(K+u),\xi)&=h(\D_{C}\!K,\xi)+\int_{S^1}\langle \alpha u,\xi\rangle dS(C,\alpha)\\
&=h(\D_{C}\!K,\xi)+\left\langle u\int_{S^1}\alpha dS(C,\alpha),\xi\right\rangle\\
&=h(\D_{C}\!K,\xi).
\end{align*}

Finally, the $\SL(W,\C)$-covariance is obtained from \eqref{eq_adjoint}. For each $g
\in\SL(W,\mathbb{C})$ we have
\begin{align*}
h(\D_C(gK),\xi)&=\int_{S^1}h(\alpha gK,\xi)dS(C,\alpha)=\int_{S^1}h(\alpha K,g^*\xi)dS(C,\alpha)\\
&=h(\D_{C}\!K,g^*\xi)=h(g\D_C\!K,\xi).
\end{align*}
It follows that $\D_C(gK)=g\D_C\!K$, hence $\D_C$ has all the required properties. 

Let us now show the uniqueness of $C$ up to translations. As the area measure $S(C,\cdot)$ is invariant under translations, we can assume that the Steiner point of $C$ is the origin. (Recall that the {\em Steiner point of a convex body $K\in\K(V)$} is defined by, see \cite[p.\! 42]{schneider_book93} 
$$s(K)=\frac{1}{\vol(B^n)}\int_{S^{n-1}}h(K,u)u du,$$
where $B^n$ denotes the unit ball in $V$.)

Let $C_{1}, C_{2}$ be convex bodies in $\C$ with $s(C_{1})=s(C_{2})=0$ and $\D_{C_{1}}= \D_{C_{2}}$, i.e. $$h(\D_{C_{1}}\!K,\xi)=h(\D_{C_{2}}\!K,\xi),\quad \forall K\in\K(W), \xi\in W^*.$$
Fix $\xi\in W^*$ and $u\in W$ such that $\xi(u)=1$. Consider $\overline{C_{1}}u\subset W$ and $\overline{C_{2}}u\subset W$. For these convex sets and $i,j\in\{1,2\}$ we have
\begin{align*}h(\D_{C_{i}}(\overline{C_{j}}u),\xi)&=\int_{S^1}h(\alpha \overline{C_{j}}u,\xi)dS(C_{i},\alpha)=\int_{S^1}h(\overline C_{j},\overline\alpha)dS(C_{i},\alpha)\\&=\int_{S^1}h(C_{j},\alpha)dS(C_{i},\alpha)=V_{2}(C_{i},C_{j}),\end{align*}
where $V_{2}$ denotes the mixed volume in $\C$. Hence, we have $V_{2}(C_{1},C_{1})=V_{2}(C_{2},C_{2})=V_{2}(C_{1},C_{2})$.

In particular, either $C_{1}$ and $C_{2}$ both have empty interior or both have non-empty interior.

Assume that $C_{1}$ and $C_{2}$ have non-empty interior. The Minkowski inequality in dimension 2 states that (see \cite[Theorem 6.2.1]{schneider_book93}) 
$$V_{2}(C_{1},C_{2})^2\geq V_{2}(C_{1},C_{1})V_{2}(C_{2},C_{2}),$$ with equality if and only if $C_{1}$ and $C_{2}$ are homothetic. Thus, we can write $C_{1}=rC_{2}+z$ with $r\in\R_{>0}$, $z\in\C$. But, from $V(C_{1},C_{1})=V(C_{2},C_{2})$ we get $r=1$ and from $s(C_{1})=s(C_{2})$, we get $z=0$. That is, $C_{1}=C_{2}$.

Assume now that $C_{1}, C_{2}$ have empty interior. Then, $C_{1}=[-z_{1},z_{1}]$ and $C_{2}=[-z_{2},z_{2}]$ with $C_{1},C_{2}\in\C$. In this case, the area measure of $C_{1}$ is given by 
$$S([-z_{1},z_{1}],\cdot)=\delta_{iz_{1}}(\cdot)+\delta_{-iz_{1}}(\cdot),$$ 
and $$h(\D_{C_{1}}\!K,\xi)=h(K,i\overline z_{1}\xi)+h(K,-i\overline z_{1}\xi).$$
Then, for every $K=[-zu,zu]$, $z\in\C$, we have
$$h(\D_{C_{1}}\!K,\xi)=h(K,i\overline z_{1}\xi)+h(K,-i\overline z_{1}\xi)=2|\!\re(i\overline z_{1}\xi(zu))|=2|\!\re(i\overline z_{1}z)|,$$
and similarly $$h(\D_{C_{1}}\!K,\xi)=h(\D_{C_{2}}\!K,\xi)=2|\!\re(i\overline z_{2}z)|.$$
It follows that $z_{1}=z_{2}$, and $C_{1}=C_{2}$.

\bigskip
Conversely, let us suppose that $Z$ is a translation invariant continuous Minkowski valuations which is $\SL(W,\C)$-covariant. We want to show that there exists some compact convex $C \subset \mathbb{C}$ with $Z=\D_C$ and $s(C)=0$.

First, we prove that $Z$ must be homogeneous of degree one. McMullen's decomposition \eqref{eq_mcmullen_dec} applied to $Z$ gives the decomposition  
\begin{displaymath}
 h(ZK,\cdot)=\sum_{k=0}^{2m}f_k(K,\cdot),
\end{displaymath}
with $f_k(K,\cdot)$ a $1$-homogeneous function. In general, 
$f_{k}$ is not subadditive as was recently proved in \cite{parapatits.wannerer}. For the 
minimal index $k_{0}$ and the maximal index $k_{1}$ with $f_{k}\neq 0$, it was proved in \cite{schneider_schuster06} that 
$f_{k_0}$ and $f_{k_1}$ are support functions.

By Lemma \ref{lemma_homog} and Remark \ref{remarkHomog} there is no non-trivial, continuous, translation invariant and $\SL(W,\mathbb{C})$-covariant Minkowski valuation $Z$ of degree $k\neq 1$, if $\dim W\geq 3$. We thus get $k_0=k_1=1$, and $Z$ is of degree 1.

\medskip
For every $\xi\in W^*$, $h(Z\cdot,\xi)$ is a real-valued valuation, which is also continuous, translation invariant and homogeneous of degree 1. Thus, by Theorem \ref{theoGoodeyWeil}, there exists a distribution $T_{\xi}$ defined on $W^*$ such that $$h(ZK,\xi)=T_{\xi}(h_{K}).$$ 

In order to derive the result for the 1-homogeneous case, we divide the proof in several steps. The first step is to show that, in our case, the distribution $T_{\xi}$ can be interpreted as a distribution on $S^1$. Then, using the $\SL(W,\C)$-covariance, we show that this distribution on $S^1$ is independent of $\xi$. The fourth step is to prove that this distribution is given by a measure defined on $S^1$. In the last two steps we find that this measure must be positive and the surface area measure of a convex set in $\C$.

\medskip
\textsc{Step 1:}
{\em Let $\xi_{0}\in W^*$. We claim that there exists a distribution $T$ on $S^1$ satisfying $(m_{\xi_{0}})_{*}T=T_{\xi_{0}}$, where $m_{\xi_{0}}:S^1\to W^*$, $m_{\xi_{0}}(\alpha)=\alpha^*\xi_{0}$, and $(m_{\xi_{0}})_{*}T(f):=T(f\circ m_{\xi_{0}})$ for every $f\in \mathcal{D}$, i.e. a continuous, 1-homogeneous function on $W^*$, smooth on $W^*\setminus\{0\}$.}

Let $E\subset W^*$ be the 1-dimensional complex subspace spanned by $\xi_{0}$.
Let $f$ be a function defined on $W^*$ such that $f|_{E}\equiv 0$. 

Let us suppose first that $f$ is the support function of $K\in\K(W)$. Then, the condition $f|_{E}=h_{K}|_{E}\equiv 0$ implies that the convex body $K$ lies in the complex subspace $F=\ker{\xi_{0}}\subset W$. 
Let $g_{\lambda}\in \GL(W,\C)$ with $g_{\lambda}^*\xi_{0}=\xi_{0}$ and $g_{\lambda}(v)=\lambda v$, $\lambda\in\R_{>0}$, for every $v\in F$. As $g_{\lambda}$ has real entries and $\det g_{\lambda}> 0$, there exist $t>0$ and $g_{0}\in\SL(W,\C)$ (with real entries) such that $g_{\lambda}=tg_{0}$. From the 1-homogeneity and the $\SL(W,\C)$-covariance of $Z$, it easily follows that $$Z(g_{\lambda}K)=g_{\lambda} ZK.$$
From the properties of $Z$ and the above equality, we get
$$h(ZK,\xi_{0})=h(ZK,g_{\lambda}^*\xi_{0})=h(Z(g_{\lambda}K),\xi_{0})=h(Z(\lambda K),\xi_{0})=\lambda h(ZK,\xi_{0}).$$
As the above equation holds for every $\lambda\in\R_{>0}$, it follows that $T_{\xi_{0}}(h_{K})=h(ZK,\xi_{0})=0$. 

Let now $f=h_{K}-h_{L}$ with $K,L\in\K(W)$ and $f|_{E}\equiv 0$, that is, $h(K,\alpha^*\xi_{0})=h(L,\alpha^*\xi_{0})$ for every $\alpha\in S^1$.
Let $g_{\lambda}\in\GL(W,\C)$ be as above. Then, $h(g_{\lambda}K,\alpha^*\xi_{0})=h(g_{\lambda}L,\alpha^*\xi_{0})$ for all $\lambda>0$. Thus,
$$\lim_{\lambda\to 0}g_{\lambda}K=\lim_{\lambda\to 0}g_{\lambda}L,$$ and from the continuity of $Z$ we get on one hand $$\lim_{\lambda\to 0}Z(g_{\lambda}K)=\lim_{\lambda\to 0}Z(g_{\lambda}L).$$

On the other hand, we have for every $\lambda\in\R_{>0}$,
$$h(ZK,\xi_{0})-h(ZL,\xi_{0})=h(Z(g_{\lambda}K),\xi_{0})-h(Z(g_{\lambda} L),\xi_{0}).$$
Taking limits on both sides we get $h(ZK,\xi_{0})=h(ZL,\xi_{0})$ and $T_{\xi_{0}}(f)=0$. 

As every function $f\in\mathcal{D}$ can be written as the difference of two support functions, that is, $f=h_{K}-h_{L}$ for some $K,L\in\K(W)$ (cf. Lemma \ref{differenceSupport}), we get $T_{\xi_{0}}(f)=0$ for every $f\in\mathcal{D}$.

Thus, we get that the value of $T_{\xi_{0}}(f)$ only depends on $f|_{E}$. We define the distribution $T$ on $S^1$ by $T(g):=T_{\xi_{0}}(\tilde g)$, where $\tilde g$ denotes an extension on $\mathcal{D}$ of $g$ (satisfying $\tilde g(\alpha^*\xi_{0})=g(\alpha)$ and hence, $T$ is well-defined). By definition of $m_{\xi_{0}}$ we have $T(f\circ m_{\xi_{0}})= T_{\xi_{0}}(f)$ since $f$ is an extension of $f\circ m_{\xi_{0}}$, and the claim follows.

\bigskip
\textsc{Step 2:} {\it Let $g\in\SL(W,\C)$ and $\xi\in W^*$. The distribution $T_{\xi}$ satisfies $T_{g^*\xi}=(g^*)_{*}T_{\xi}$, where $(g^*)_{*}T_{\xi}(f)=T_{\xi}(f\circ g^*)$ for every $f\in\mathcal{D}$.}

We first prove the equality for a support function $h_{K}$, $K\in\K(W)$. Using property \eqref{eq_adjoint} of support functions, and that $Z$ is an $\SL(W,\C)$-covariant valuation, we get 
\begin{align*}(g^*)_{*}T_{\xi}(h_{K})&=T_{\xi}(h_{K}\circ g^*)=T_{\xi}(h_{gK})\\&=h(Z(gK),\xi)=h(ZK,g^*\xi)=T_{g^*\xi}(h_{K}).\end{align*}

The general case follows by linearity and Lemma \ref{differenceSupport}.

\bigskip
\textsc{Step 3:}
{\it The distribution $T$ on $S^1$ given in Step 1 satisfies $(m_{\xi})_{*}T=T_{\xi}$, for every $\xi\in W^*$.}

Let $\xi_{0}\in W^*$ as in Step 1 and $\xi\in W^*$. There exists $g\in \SL(W,\C)$ such that $g^*\xi_{0}=\xi$. Using Steps 2 and 1, it follows that
\begin{align*}T_{\xi}=T_{g^*\xi_{0}}=(g^*)_{*}T_{\xi_{0}}=(g^*)_{*}(m_{\xi_{0}})_{*}T=(g^*\circ m_{\xi_{0}})_{*}T=(m_{\xi})_{*}T.
\end{align*}

\bigskip
\textsc{Step 4:}
{\it The distribution $T$ defined in Step 1 is given by a signed measure $\mu$. That is, $$(m_{\xi})_{*}T(f)=\int_{S^1}f(\alpha^*\xi)d\mu(\alpha).$$}

Schneider obtained in \cite{schneider74} a classification of continuous, Min\-kows\-ki valuations $\Phi$ on a 2-dimensional vector space $V$ which satisfy $\Phi b = b \Phi$ for every $b$ in $\SO(V,\R)$ or $b$ a translation in $V$. The general expression for such a $\Phi$ is 
$$h(\Phi(K),\alpha)=\int_{0}^{2\pi}h(K-s(K),u(\alpha+\beta))d\nu(\beta)+\langle u(\alpha),s(K)\rangle,$$
where $s(K)$ denotes the Steiner point of $K$, $u(\alpha)=\cos(\alpha)e^{1}+\sin(\alpha)e^{2}$ with $\{e^{1},e^{2}\}$ a given basis on $V^*$, and $\langle\,,\rangle$ the pairing of $V$ with its dual space. The signed measure $\nu$ is unique up to a linear measure $l$ defined by $\omega\mapsto\int_{\omega}\langle u(\alpha),a\rangle d\alpha$ with $a\in V$ a constant vector. 

Schneider's proof can be easily adapted to our situation. Indeed, let $\xi\in W^*$, $\{e^1,\dots,e^m\}$ a basis of $W^*$ with $e^1=\xi$ and $\{e_{1},\dots,e_{m}\}$ the basis of $W$ dual to $\{e^1,\dots,e^m\}$. Define $E$ as the 1-dimensional complex space in $W$ spanned by $e_{1}$. We can identify $E^*$ with $\mathrm{span}_{\C}(\xi)\subset W^*$, and write $\alpha\in E^*$ as a multiple of $\xi$, which we denote again by $\alpha$.
 
Let $\phi:E\to E$ be the restriction of $Z$ to $E$, i.e. $h(\phi(K),\alpha)=h(ZK,\alpha\xi)$ with $K\in\K(E)$.

The operator $\phi$ inherits the properties of $Z$, that is, it is a continuous,  translation invariant Minkowski valuation which is covariant with respect to $\SL(E,\C)\cong S^1$. Now, $\phi$ satisfies all the hypothesis of Schneider's result except the covariance with respect to translations (i.e. $\Phi(K+t)=\Phi(K)+t$), which is replaced by the invariance (i.e. $\phi(K+t)=\phi(K)$). However, the first step in Schneider's proof is to construct a translation invariant valuation from the translation covariant one via $\Phi-s$ (where $s$ denotes the Steiner point), and the same argument can be used in our situation.

\bigskip
\textsc{Step 5:}
{\it The measure $\mu$ is positive.}

By assumption, the function $F(\xi)=h(ZK,\xi)$ must be convex. Thus, the second differential of $F$ at each $\xi$ must be a positive semi-definite bilinear form (cf. \cite[p.\! 108 and Theorem 1.5.10]{schneider_book93}). Note that, by the representation of $h(ZK,\xi)$ obtained in the previous step, if $h(K,\cdot)$ is a smooth function, then $h(ZK,\cdot)$ is also smooth.

For simplicity, we fix a scalar product on $W$. Fix $\xi\in W$ and let $S^1\cdot\xi\subset W$ be a circle contained in the complex line spanned by $\xi$. Let $\epsilon>0$ and $\mathring{B}$ the open $2m$-dimensional ball in $W$. We have that $S^1\cdot\xi\times \epsilon \mathring{B}$ is a neighborhood of $S^{1}\cdot \xi$ in $W$. Geometrically, it can be interpreted as the open tube of radius $\epsilon$ along $S^1\cdot\xi$.

Let $f:S^1\to\R_{>0}$ be a positive smooth function. 
Attach to each point $p_{\theta}:=e^{i\theta}\xi\in S^1\cdot\xi$ a spherical cap of a $(2m-2)$-dimensional sphere $S_{\theta}$ with radius $f(\theta)$, center on the segment $[0, p_{\theta}]$ and tangent plane at $p_{\theta}$ orthogonal to $\{\xi, i\xi\}$. 

For $\epsilon$ small enough the intersection between $S^1\cdot\xi\times\epsilon \mathring{B}$ and the set described in the previous paragraph is a smooth hypersurface. Denote it by $K_{\epsilon}$.

The principal curvatures of $K_{\epsilon}$ at $p_{\theta}$ are $f^{-1}(\theta)>0$ (the inverse of the radius of the attached $(2m-2)$-dimensional sphere at $p_{\theta}$) with multiplicity $2m-2$, and 1, corresponding to the principal direction $ie^{i\theta}\xi$, tangent to $S^{1}\cdot\xi$ at $p_{\theta}$. Since the function $f$ is smooth, strictly positive, and defines the hypersurface in a neighborhood of $S^1\cdot\xi$, all principal curvatures at any point of $K_{\epsilon}$ are strictly positive, provided $\epsilon$ is small enough.

Set $K:=\mathrm{conv} \overline{K_{\epsilon}}$, where $\mathrm{conv}\overline{ K_{\epsilon}}$ denotes the convex hull of $\overline{K_{\epsilon}}$, the closure of $K_{\epsilon}$, at $W$. As $K_{\epsilon}$ is a convex hypersurface, we have $\partial K\cap(S^1\cdot\xi\times\epsilon B)=K_{\epsilon}\cap(S^1\cdot\xi\times\epsilon B)$, and $K$ is smooth in a neighborhood of $S^{1}\cdot \xi$.
The second differential of $h(ZK,\cdot)$ at $\xi$ must be positive semi-definite (cf. \cite[Theorem 1.5.10]{schneider_book93}), that is,
$$(d^2h(ZK,\xi))(a,a)=\int_{S^1}(d^2h(\alpha K,\xi))(a,a)d\mu(\alpha)\geq 0,\,\forall a\in W.$$ 

The eigenvalues of the second differential of the support function of a convex body $K\subset W$ in a direction $\xi\in W$ are the radii of curvature of $K$ at the corresponding supporting point with eigenvector the principal directions at this supporting point (cf.\! \cite[Corollary 2.5.2]{schneider_book93}). By the construction of $K$, the support point of $\overline\alpha K$ in direction $\xi$ is $p_{-\theta}$, with $e^{-i\theta}=\overline\alpha$.

Take $a=b=u$ with $u$ a principal direction of $K$ at $p_{0}=\xi$ different from $i\xi$. Then, $u$ is also a principal direction of $K$ at every point $p_{\theta}$, with principal radius of curvature $f(e^{i\theta})$, the radius of the attached sphere at $p_{\theta}$. 

Hence, 
$$(d^2h(\alpha K,\xi))(u,u)=f(\overline\alpha).$$ 
Therefore, for every strictly positive smooth function on $S^1$ we get $$\int_{S^1}f(\overline\alpha)d\mu(\alpha)\geq 0,$$
and from the density of smooth positive functions on the space of continuous positive functions, we get that the measure $\mu$ is positive. 

\bigskip
\textsc{Step 6:}
{\it The measure $\mu$ is the surface area measure for some convex body in $\C$.}

Using the translation invariance of $Z$, we have that $\mu$ satisfies condition \eqref{zero} of Minkowski's existence theorem (see Theorem \ref{minkowskiExistence}). 

If $\mu$ is not concentrated on two antipodal points of $S^1$, then from Theorem \ref{minkowskiExistence} we get the existence of a 2-dimensional convex body $C\subset \C$ with $dS(C,\cdot)=\mu(\cdot)$.

Otherwise, if $\mu$ is concentrated on two antipodals points $\pm \alpha$ of $S^1$, then $\mu$ coincides with the surface area of a centered interval with normal vector given by the direction $\alpha$. Thus, $\mu$ is the surface area of a (1-dimensional) convex body.  
\endproof

\section{The case $\dim W=2$}\label{dim2cov}
In order to have a complete classification of continuous, translation invariant, $\SL(W,\C)$-covariant or $\SL(W,\C)$-contravariant Minkowski valuations with a fixed homogeneity degree, it just remains to study the case of $\SL(W,\C)$-covariant valuations of degree 3 and $\SL(W,\C)$-contravariant valuations of degree 1 in a complex 2-dimensional space $W$.

Fix a basis of $W$ and consider the determinant map 
\begin{equation}\label{detId}\func{\det}{W\times W}{\C}{(u,v)}{\det(u,v).}
\end{equation}
This map induces an identification $\Phi$ between $W$ and its dual space $W^*$, which satisfies $\Phi(gu)=(\det g) g^{-*}\Phi(u)$, for every $g\in\GL(W,\C)$, $u\in W$.

Then, every $\SL(W,\C)$-contravariant (resp. covariant) Minkowksi valuation $Z$ of degree $k$ is in correspondence with an $\SL(W,\C)$-covariant (resp. contravariant) Minkowski valuation $\Phi^{-1}\circ Z$ (resp. $\Phi\circ Z$) also of degree $k$. Thus, the following classification results follow directly from Theorem 1.1 and Theorem 1.

\begin{Proposition}Let $\dim_{\C}W=2$ and $Z:\K(W)\to\K(W^*)$ a continuous, translation invariant and $\SL(W,\C)$-contravariant Minkowski valuation of degree $k$. If $k\neq 1,3$, then $Z\equiv \{0\}$. If $k=3$, then $Z$ is of the form \eqref{complexProjection}. If $k=1$, then there exists a convex body $C\subset \C$ such that 
\begin{displaymath}h(ZK,w)=\int_{S^1}h(\det(K,w),\alpha)dS(C,\alpha),\quad K\in\K(W),\,w\in W,\end{displaymath}
where $\det(K,w):=\{\det(k,w)\,|\,k\in K\}$. Moreover, $C$ is unique up to translations. 
\end{Proposition}

\begin{Proposition}Let $\dim_{\C}W=2$ and $Z:\K(W)\to\K(W)$ a continuous, translation invariant and $\SL(W,\C)$-covariant Minkowski valuation of degree $k$. If $k\neq 1,3$, then $Z\equiv \{0\}$. If $k=1$, then $Z$ is of the form \eqref{eq_def_pi_c}. If $k=3$, then there exists a convex body $C\subset \C$ such that 
\begin{displaymath}\label{cova1new2}h(ZK,\xi)=V(K,K,K,C\cdot w),\end{displaymath}
where $w\in W$ is the corresponding vector to $\xi$ given by the identification $\Phi^{-1}$ between $W^*$ and $W$ determined by \eqref{detId}.
Moreover, $C$ is unique up to translations.
\end{Proposition}

\bibliographystyle{plain}
\bibliography{biblio}
\end{document}